 \documentstyle[graphicx,amscd,amssymb,verbatim,12pt,righttag,color]{amsart}
 \newtheorem{theorem}{Theorem}[section]
 \newtheorem{Def}[theorem]{Definition}
 \newtheorem{Prop}[theorem]{Proposition}
 \newtheorem{Lem}[theorem]{Lemma}
 
 \newtheorem{Rem}[theorem]{Remark}
 \newtheorem{Example}[theorem]{Example}

\newcommand{\R}{{\mathbb R}}
\newcommand{\Z}{{\mathbb Z}}

\newcommand{\N}{{\mathbb N}}

\newcommand{\C}{{\mathcal C}}
\newcommand{\D}{{\mathcal D}}

\newcommand{\M}{{\mathcal M}}

\newcommand{\var}{{\vartheta}}

\newcommand{\eproof}{\hfill$\square$}%\rule{2.2mm}{3.0mm}}

\newcommand{\proof}{\noindent {\bf Proof.~~}}

 \numberwithin{equation}{section}
 \renewcommand{\rm}{\normalshape}
\begin{document}

\title  {Spectral property of Cantor measures  with consecutive digits}

%\date{29 November 2011}
\author{Xin-Rong Dai}
\address{School of Mathematics and Computational Science, Sun Yat-Sen University, Guangzhou, 510275,  P. R. China}
\email{daixr@@mail.sysu.edu.cn}
\author{Xing-Gang He}
\address{College of Mathematics and Statistics, Central China Normal
University, Wuhan 430079, P. R. China}
 \email{xingganghe@@163.com}
\author{Chun-Kit Lai}
\address{Department of Mathematics and Statistics, McMaster University,
Hamilton, Ontario, L8S 4K1, Canada}
\email{cklai@@math.mcmaster.ca}

\thanks{The research is supported in part by the HKRGC Grant and the Focused Investment Scheme of CUHK; The second author is also supported by
the National Natural Science Foundation of China 11271148.}

%\date{\today}
%\keywords { Blocking, cyclotomic polynomials, kernel polynomials, prime, product-forms, self-similar tiles, spectra, tile digit sets, tree. }
%\thanks{ The research is supported in part by the HKRGC Grant and the Focused Investment Scheme of CUHK. }
\maketitle

\begin{abstract}
We consider  equally-weighted Cantor measures $\mu_{q,b}$
arising from iterated function systems of the form ${b^{-1}(x+i)}$,
$i=0,1,\cdots,q-1$, where $q<b$. We classify the $(q,b)$ so that
they have infinitely many mutually orthogonal exponentials in $L^2(\mu_{q,b})$. In particular, if $q$ divides $b$, the
measures have a complete orthogonal exponential system
and hence spectral measures. Improving the construction in
\cite{[DHS]}, we characterize all the maximal orthogonal sets $\Lambda$ when
$q$ divides $b$ via a maximal mapping on the $q-$adic tree in which all elements in $\Lambda$
are represented uniquely in finite $b-$adic expansions and we can separate
 the maximal orthogonal sets into two
 types: regular and irregular sets. For a regular maximal orthogonal set, we show that its completeness
  in $L^2(\mu_{q,b})$ is crucially determined
 by the certain growth rate of non-zero digits in the tail of the $b-$adic expansions of the elements. Furthermore, we exhibit complete orthogonal
 exponentials with zero Beurling dimensions. These examples show that the technical condition in Theorem 3.5 of \cite{[DHSW]} cannot be removed.
For an irregular maximal orthogonal set, we show that under some condition, its completeness is equivalent to that of the corresponding regularized mapping.
\end{abstract}

\tableofcontents
%\bigskip

\section{Introduction}

Let $\mu$ be a compactly supported Borel probability measure on
${\Bbb R}^d$. We say that $\mu$ is a {\it spectral measure} if there
exists a countable set $\Lambda\subset {\Bbb R}^d$ so that
$E(\Lambda): = \{e^{2\pi i \langle\lambda,x\rangle}:
\lambda\in\Lambda\}$ is an orthonormal basis for $L^2(\mu)$. In this
case, $\Lambda$ is called a {\it spectrum} of $\mu$. If
$\chi_{\Omega}dx$ is a spectral measure, then we say that $\Omega$
is a {\it spectral set}. The study of spectral measures
was first initiated from B. Fuglede
in 1974 \cite{[Fu]}, when he considered a functional analytic
problem of extending some commuting partial differential operators
to some dense subspace of $L^2$ functions. In his first attempt,
Fuglede proved that any fundamental domains given by  a discrete
lattice are spectral sets with its dual lattice as its spectrum.  On
the other hand, he also proved that triangles and circles on ${\Bbb
R}^2$ are not spectral sets, while some examples (e.g.
$[0,1]\cup[2,3]$) that are not fundamental domains can still be
spectral. From the examples and the relation
between Fourier series and translation operators, he proposed a
reasonable conjecture on spectral sets: {\it $\Omega\subset{\Bbb
R}^s$ is a spectral set if and only if $\Omega$ is a translational
tile.} This conjecture baffled  experts for 30
years until 2004, Tao \cite{[T]} gave the first counterexamples on ${\Bbb R}^d$,
$d\geq5$. The examples were modified
later so that the conjecture are false in both directions on ${\Bbb
R}^d$, $d\geq3$ \cite{[KM1]}, \cite{[KM2]}. It remains open in
dimension 1 and 2. Despite the counterexamples, the
exact relationship between spectral measures and tiling
is still mysterious.

\medskip

The problem of spectral measures is as exciting when we consider fractal measures. Jorgensen and Pedersen \cite{[JP]} showed that the standard
 Cantor measures are spectral measures if the contraction is $\frac{1}{2n}$, while there are at most two orthogonal
 exponentials when the contraction
 is $\frac{1}{2n+1}$. Following this discovery, more spectral self-similar/self-affine measures were also found (\cite{[LaW]}, \cite{[DJ]}
 et. al.).
 The construction of these spectral self-similar measures is based on the existence of the {\it compatible pairs
  (known also as Hadamard triples)}.
 It is still unknown whether all such spectral measures are obtained from compatible pairs. Having an exponential basis, the series convergence
 problem was also studied by Strichartz. It is surprising that the ordinary Fourier
  series of continuous functions converge uniformly for
 standard Cantor measures \cite{[Str]}. By now there are considerable amount of literatures studying spectral measures and other generalized types
 of Fourier expansions like the Fourier frames and Riesz bases (\cite{[DHJ]}, \cite{[DHSW]}, \cite{[DL]}, \cite{[HLL]}, \cite{[IP]}, \cite{[JKS]}
 \cite{[La]}, \cite{[LaW]}, \cite{[Lai]}, \cite{[Li1]}, \cite{[Li2]}, and the references therein).

\medskip

In \cite{[HuL]}, Hu and Lau  made a start in studying the spectral
properties of Bernoulli convolutions, the simplest class of
self-similar measures.
 They classified the contraction ratios with infinitely many orthogonal exponentials. It was recently shown by Dai  that the only spectral
  Bernoulli convolutions are of contraction ratio $\frac{1}{2n}$ \cite{[D]}. In this paper, we study another general class of
 Cantor
  measures on ${\Bbb R}^1$. Let $b>2$ be an integer and $q<b$ be another positive integer. We consider the iterated function system(IFS) with maps
$$
f_i(x) = b^{-1}(x+i), \ i = 0,1,...,q-1.
$$
The IFS arises a natural {\it self-similar measure} $\mu = \mu_{q,b}$ satisfying
\begin{equation}\label{eq1.1}
\mu (E) = \sum_{i=0}^{q-1}\frac{1}{q}\mu( f_j^{-1}(E))
\end{equation}
for all Borel sets $E$. Note that we only need to consider equal weight since non-equally weighted self-similar measures here cannot have any spectrum
by Theorem 1.5 in \cite{[DL]}. It is also clear that if $q=2$, $\mu$ becomes the standard Cantor measure of $b^{-1}$ contraction. For this class
 of self-similar measures, we find surprisingly that the spectral properties depend heavily on the number theoretic relationship between $q$ and $b$.
  Our first result is to show that {\it $\mu = \mu_{q,b}$ has infinitely many orthogonal exponentials if and only if $q$ and $b$ is not
  relatively prime. If moreover, $q$ divides $b$, the resulting measure will be a spectral measure} (Theorem \ref{th1.1}). However,
  when $q$ does not divide $b$ and they are not relatively prime (e.g. $q=4, \ b=6$), variety of cases may occur and we are not sure
  whether there are spectral measures in these classes (see Remark \ref{rem2.1}).

\medskip

We then focus on the case when $b = qr$ in which we aim at giving a
detailed classification of its spectra. The classification of
spectra,  for a given spectral measure, was first studied by
Lagarias, Reeds and Wang \cite{[LRW]}. They considered the spectra
of $L^2([0,1)^d)$ (more generally fundamental domains of some
lattices) and they showed that the spectra of $L^2([0,1)^d)$ are
exactly all the tiling sets of $[0,1)^d$. If $d=1$, the way of
tiling  $[0,1)$ is rather rigid, and it is easy to see that the only
spectrum (respectively the tiling set) is the translates of the
integer lattice ${\Bbb Z}$.

\medskip

Such kind of rigidity breaks down even on ${\Bbb R}^1$ if we turn to
fractal measures. The first attempt of the classification of its spectra was due
 to \cite{[DHS]}, Dutkay, Han and Sun decomposed the maximal orthogonal sets of one-fourth Cantor measure using $4$-adic expansion with
 digits $\{0,1,2,3\}$ and put them into a labeling of the binary tree. The maximal orthogonal sets will then be obtained by reading all
 the infinite paths with digits ending eventually in $0$ (for positive elements) or $3$ (for negative elements). They also gave some
 sufficient conditions on the digits for a maximal orthogonal set to be a spectrum. Nonetheless, the condition is not easy to verify.

\medskip

 Turning to our self-similar measures with consecutive digits where the one-fourth Cantor measure is a special case, we will classify all the maximal
  orthogonal sets using mappings on the standard $q-$adic tree called {\it maximal mappings} (Theorem \ref{th1.6}). This construction improves the
  tree labeling method in \cite{[DHS]} in two ways.

\begin{enumerate}
\item We will choose the digit system to be $\{-1,0,1,\cdots, b-2\}$ instead of $\{0,1,\cdots,$ $b-1\}$. By doing so, all integers
(both positive and negative) can be expanded into $b-$adic
expansions terminating at $0$ (Lemma \ref{lem1.1}).
\item We impose restrictions on  our  labeling position on the tree so that together with (1), all the elements
in a maximal orthogonal set can be extracted by reading some
specific paths in the tree. These
paths are collected in a countable set $\Gamma_q$  defined in
(\ref{gamma}).
\end{enumerate}

Having such a new tree structure of a maximal orthogonal set, we discover there are
two possibilities for the maximal sets depending on
whether all the paths in $\Gamma_q$ are corresponding to some
elements in the maximal orthogonal sets (i.e. the values assigned
are eventually 0). If it happens that all the paths in $\Gamma_q$
behave nicely as said, we call such maximal orthogonal sets {\it
regular}. It turns out that regular
 sets cover most of the interesting cases and we can give regular sets a natural ordering $\{\lambda_n: n=0,1,2,\cdots\}$.
 If the standard $q-$adic expansion  of $n$  has length $k$, we define $N_n^{\ast}$ to be the number of non-zero digits in the $b-$adic
 expansion using $\{-1,0,\cdots,b-2\}$ of $\lambda_n$ after $k$ . $N_n^{\ast}$ is our crucial factor in determining whether the set is a spectrum.
 We show that {\it if $N_n^{\ast}$ grows slowly enough or even uniformly bounded, the set will be a spectrum, while if  $N_n^{\ast}$ grows too fast,
 say it is of polynomial rates, then the maximal orthogonal sets will not be a spectrum} (Theorem \ref{th1.7}).

\medskip

In \cite{[DHSW]}, Dutkay {\it et al } tried to generalize the
classical results of Landau \cite{[Lan]} about the Beurling density
on Fourier frames to fractal settings. They defined the concept of
{\it Beurling dimensions} for a discrete set and showed that all
Bessel sequences for an IFS of similitudes with no overlap condition
must have Beurling dimension not greater than its Hausdorff
dimension of the attractor. Under technical assumption on the frame
spectra, they showed the above two dimensions coincide. They
conjectured that the assumption can be removed. However, as we see
that $N_n^{\ast}$ counts the number of non-zero digits only, we can
freely add $qb^m$ for any $m>0$ on the tree of the canonical
 spectrum. These additional terms
 push the $\lambda_n$'s  as far away
 from each other as wanted and we therefore show that  {\it there exists  spectrum of zero Beurling dimension} (Theorem \ref{th1.9}).

\medskip

For the organization of the paper, we present our set-up and main
results in Section 2. In Section 3, we discuss the maximal
orthogonal sets of $\mu_{q,b}$ and classify all maximal orthogonal
sets via the maximal mapping on the $q-$adic when $q$ divides $b$. In Section 4, we discuss the
regular
 spectra and prove the growth rate criteria. Moreover, the examples of the spectra with zero Beurling dimensions will be given. In Section 5,
  we give a study on the irregular spectra.

\bigskip

\section{Setup and main results}

  Let $\Lambda$ be a countable set in ${\Bbb R}$ and denote
 $E(\Lambda)=\{e_\lambda: \lambda\in\Lambda\}$ where
 $e_\lambda(x)=e^{2\pi i\lambda x}$. We say that $\Lambda$ is a {\it maximal
 orthogonal set} ({\it spectrum}) if $E(\Lambda)$ is a maximal orthogonal
 set  (an orthonormal basis) for $L^2(\mu)$. Here $E(\Lambda)$ is a maximal orthogonal
 set of exponentials means that it is a mutually orthogonal set in $L^2(\mu)$ such that if $\alpha\not\in\Lambda$, $e_{\alpha}$ is not orthogonal
 to some $e_{\lambda}$, $\lambda\in\Lambda$. If $L^2(\mu)$ admits a spectrum, then $\mu$ is called a {\it spectral measure}. Given a measure $\mu$,
 the Fourier transform is defined to be
 $$
 \widehat{\mu}(\xi) = \int e^{2\pi i \xi x}d\mu(x).
 $$
It is easy to see that $E(\Lambda)$ is an orthogonal set if and only if
$$
(\Lambda-\Lambda)\setminus\{0\}\subset {\mathcal Z}(\widehat{\mu}) := \{\xi\in{\Bbb R}: \widehat{\mu}(\xi)=0\}.
$$
We call such $\Lambda$ a {\it bi-zero set} of $\mu$.
 For $\mu = \mu_{q,b}$, we can calculate its Fourier transform.
\begin{equation}\label{eq1.2}
 \widehat{\mu}(\xi) = \prod_{j=1}^{\infty}\left[\frac{1}{q}(1+e^{2\pi i b^{-j}\xi}+...+e^{2\pi i b^{-j}(q-1)\xi})\right].
 \end{equation}
Denote
\begin{equation}\label{eq1.2+}
m(\xi) = \frac{1}{q}\left(1+e^{2\pi i \xi}+\cdots+ e^{2\pi i (q-1)\xi}\right)
\end{equation}
and thus $|m(\xi)|= |\frac{\sin q\pi \xi}{q\sin \pi \xi}|$. The zero set of $m$ is
$$
{\mathcal Z}(m) = \left\{\frac{a}{q}: q\nmid a, a\in\Z\right\},
$$
where $q\nmid a$ means $q$ does not divide $a$. We can then write $
\widehat{\mu}(\xi) = \prod_{j=1}^{\infty}m(b^{-j}\xi)$,
so that the  zero set of $\widehat{\mu}$ is given by
 \begin{equation}\label{eq1.3}
{\mathcal Z}(\widehat{\mu})   = \left\{\frac{b^{n}}{q} a:  n\geq 1,
\ q\nmid a \right\}=r\{b^na: n\ge 0, q\nmid a\},
 \end{equation}
where $r=b/q$.
 \medskip

We have the following theorem classifying which $\mu_{q,b}$ possess infinitely many orthogonal exponentials. It is also the starting point of our paper.

 \begin{theorem}\label{th1.1}
 $\mu = \mu_{q,b}$ has infinitely many orthogonal exponentials if and only if
the greatest common divisor between $q$ and $b$ is  greater than
$1$. If $q$ divides $b$, then $\mu_{q,b}$ is a spectral measure.
\end{theorem}

\bigskip

 We wish to give a classification on the spectra and the maximal orthogonal sets whenever they exist. To do this, it is convenient to introduce
 some multiindex notations:  Denote $\Sigma_q = \{0,\cdots, q-1\}$,
 $\Sigma_q^0=\{\var\}$ and
  $\Sigma_q^n= \underbrace{\Sigma_q\times\cdots\times \Sigma_q}_n$.
 Let $\Sigma_q^{\ast} = \bigcup_{n=0}^{\infty}\Sigma_q^n$ be the set of all finite words and let $\Sigma_q^{\infty} =\Sigma_q\times
 \Sigma_q\times\cdots $ be the set of all infinite words. Given $\sigma=\sigma_1\sigma_2\cdots \in \Sigma^{\infty}\cup\Sigma^{\ast}$,
 we define $\var\sigma=\sigma$, $\sigma|_{k} =\sigma_1\cdots\sigma_k$ for $k\ge 0$ where $\sigma|_0=\var$ for any $\sigma$ and adopt the
 notation ${0}^{\infty} = 000\cdots$, $0^k =\underbrace{0\cdots0}_k$ and $\sigma\sigma'$ is the concatenation of $\sigma $ and $\sigma'$.
 We start with a definition.

\begin{Def}
Let $\Sigma_q^*$ be all the finite words defined as above. We say it
is a {\it $q-$adic tree} if we set naturally the root is $\var$, all
the $k$-th level nodes are $\Sigma_q^k$ for $k\ge 1$ and all the
offsprings of $\sigma\in \Sigma_q^*$ are $\sigma i$ for $i=0,
1,\ldots, q-1$.
\end{Def}

Let $\tau$ be a map from $\Sigma_q^*$ to real numbers.
Then the image of $\tau$
defines a {\it $q$-adic tree labeling}. Define $\Gamma_q$
\begin{equation}\label{gamma}
\Gamma_q: = \{\sigma{0}^{\infty}: \  \sigma =\sigma_1\cdots\sigma_k\in\Sigma_q^*, \ \sigma_k\neq 0\}.
\end{equation}
$\Gamma_q$ will play a special role in our construction.

Suppose that for some word $\sigma= \sigma'0^{\infty}\in\Gamma_q$,
$\tau(\sigma|_k)=0$ for all $k$ sufficiently large, we say that
$\tau$ is {\it regular on $\sigma$}, otherwise {\it irregular}. Let
$b$ be another integer, if $\tau$ is regular on some $\sigma\in
\Gamma_q$, we define the projection $\Pi^\tau_b$ from $\Gamma_q$ to
$\R$ as
\begin{eqnarray}\label{formally}
\Pi^\tau_b(\sigma)=\sum_{k=1}^\infty \tau(\sigma|_k)b^{k-1}.
\end{eqnarray}
The above sum is finite since $\tau(\sigma|_k) =0$ for sufficiently
large $k$. If $\tau$ is regular on any $\sigma$ in $\Gamma_q$, we
say that $\tau$ is a {\it regular mapping}.

\medskip

\begin{Example}\label{example1.1} Suppose $b=q$,  let $\C=\{c_0=0, c_1, \ldots, c_{b-1}\}$ be a residue system mod $b$ where $c_i\equiv i$ (mod $b$).
 Define $\tau(\var)=0$ and $\tau(\sigma)=c_{\sigma_k}$ if $\sigma=\sigma_1\cdots \sigma_k\in\Sigma_q^k\subset\Sigma_q^{\ast}$. Then it is easy to see that $\tau$
 is regular on any $\sigma\in\Gamma_q$ and hence it is regular. Moreover,
\begin{equation}\label{first}
\Pi^\tau_b(\Gamma_b)\subseteq \Z.
\end{equation}
When $\C=\{0, 1, \ldots, b-1\}$,  then the mapping $\Pi^\tau_b$ is
 a bijection from $\Gamma_b$ onto $\N\cup\{0\}$.
\end{Example}

\medskip

In [DHS], putting their setup in our language, they classified
maximal orthogonal sets of standard one-fourth Cantor measure via
the mapping $\tau$ from $\Sigma_2^{\ast}$ to $\{0,1,2,3\}$. However,
some maximal orthogonal sets may have negative elements in which
those elements cannot be expressed finitely in $4$-adic expansions
using digits $\{0,1,2,3\}$. In our classification, we will choose
the digit system to be
 ${\C} = \{-1,0,1\cdots,b-2\}$ in which we can expand any integers uniquely by finite $b-$adic expansion. We have the following simple but
 important lemma.

\begin{Lem}\label{lem1.1}
Let $\C=\{-1, 0, 1, \ldots, b-2\}$ with integer $b\geq3$ and let
$\tau$ be the map defined in Example \ref{example1.1}. Then
$\Pi^\tau_b$ is
 a bijection between $\Gamma_b$ and ${\Bbb Z}$.
\end{Lem}

\proof For any $n\in \Z$ and $|n|<b$, it is easy to see that there exists unique $\sigma\in\Gamma_b$ such that
 $n=\Pi^\tau_b(\sigma)$. For example, $n=b-1$, then $n=\Pi^\tau_b(\sigma_1\sigma_2)$ where $\sigma_1=b-1$ and $\sigma_2=1$.
When $|n|\ge b$, then $n$ can be decomposed uniquely as $n=\ell b+c$ where $c\in\C$.  We note that $|\ell|=|\frac{n-c}b|\leq\frac{|n|+b-2}b<|n|$.
 If $|\ell|<b$, we are done. Otherwise, we further decompose $\ell$ in a similar way and after finite number of steps, $|\ell|<b$. The expansion is
 unique since each decomposition is unique. \eproof

\medskip

We now define a $q-$adic tree labeling  which  corresponds to a maximal orthogonal
set for $\mu_{q,b}$ when $b=qr$. We
observe that for $b=qr$, we can decompose ${\mathcal C} =
\{-1,0\cdots,b-2\}$ in $q$ disjoint classes according to the
remainders after being divided by $q$: ${\mathcal C} =
\bigcup_{i=0}^{q-1}{\mathcal C}_i$ where
$$
{\mathcal C}_i = (i+q{\Bbb Z})\cap {\mathcal C}.
$$

\medskip

\begin{Def}\label{def1.5}
{\rm Let $\Sigma_q^*$ be a $q-$adic tree and $b = qr$, we say that $\tau$ is a {\it maximal mapping} if it is a map
$\tau = \tau_{q,b}: \Sigma_q^*\rightarrow \{-1,0,...,b-2\}$ that satisfies }

{\rm (i) $\tau(\vartheta) =\tau(0^n) =0$  for all $n\geq1$.}

{\rm (ii) For all $k\geq1$}, $\tau (\sigma_1\cdots \sigma_k) \in {\mathcal C}_{\sigma_k}.$

{\rm (iii)  For any word $\sigma\in \Sigma_q^{\ast}$, there exists
$\sigma'$ such that  $\tau$ is regular on
$\sigma\sigma'0^{\infty}\in\Gamma_q$.}
\end{Def}

 We call a tree mapping a {\it regular mapping} if it
satisfies (i) and (ii) in above and is regular on any word in
$\Gamma_q$. Clearly, regular mappings are maximal.

\bigskip

 Given a maximal mapping $\tau$,   the following
 sets will be of our main study in this paper.
\begin{equation}\label{eq1.4}
\Lambda(\tau):=\{\Pi^\tau_b(\sigma): \sigma\in\Gamma_q ,\ \mbox{
$\tau$ is regular on $\sigma$}\}.
\end{equation}

\begin{figure}[h]
\centerline{\includegraphics[width=15cm,height=7.5cm]{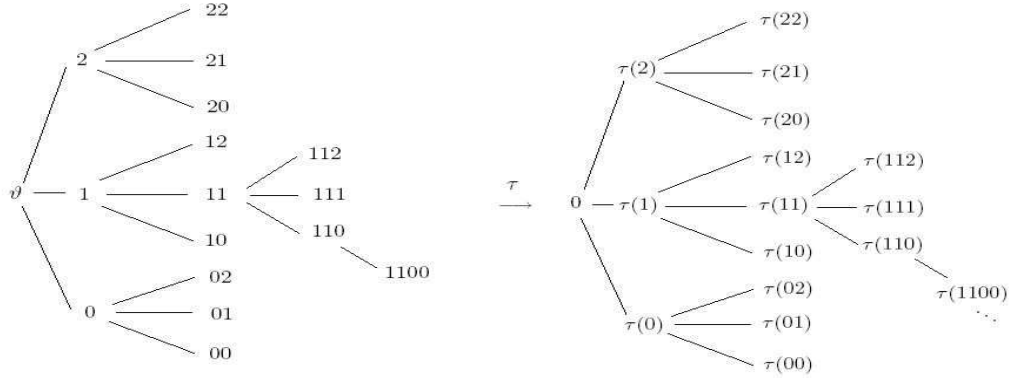}}
\caption {\small{ An illustration of the $3-$adic tree and the
associated mapping $\tau$.}}
\end{figure}

\medskip

From now on,  we will assume that $b=qr$, $\C=\{-1, 0, 1,\ldots,
b-2\}$ and $0\in \Lambda$. The main results are as follows and this
is also the reason why $\tau$ is called a maximal mapping.

\begin{theorem}\label{th1.6}
$\Lambda$ is a maximal orthogonal set of $L^2(\mu_{q, b})$ if and
only if there exists a maximal mapping $\tau$ such that $\Lambda = r
\Lambda(\tau)$, where $b=qr$.
\end{theorem}

For the proof, (i) in Definition \ref{def1.5} is to ensure $0\in\Lambda$. (ii) is to make sure the mutually orthogonality and (iii) is for the
maximal orthogonality.

\bigskip
If $\Lambda$ is a spectrum of $L^2(\mu)$, we call the associated
maximal mapping $\tau$ a {\it spectral mapping}. We will restrict
our attention to regular mappings (i.e. for all $\sigma\in\Gamma_q$,
$\tau$ is regular on $\sigma$). In this case, $\Lambda(\tau) =
\{\Pi^\tau_b(\sigma): \sigma\in\Gamma_q \}$.
 The advantage of considering regular mappings is that we can give a natural ordering of the maximal orthogonal set $\Lambda(\tau)$.
 The ordering goes as follows: Given any $n\in{\Bbb N}$, we can expand it into the unique finite $q-$adic expansion,
\begin{equation}\label{eq1.5}
n= \sum_{j=1}^{k}\sigma_jq^{j-1}, \
\sigma_j\in\{0,\cdots,q-1\},\qquad \sigma_k\ne 0.
\end{equation}
 In this way  $n$ is uniquely corresponding to one word $\sigma=\sigma_1\cdots\sigma_k$, which is called the {\it $q$-adic expansion of $n$}.
  For a regular  mapping $\tau$, there is a natural ordering of the maximal orthogonal set $\Lambda(\tau)$:  $\lambda_{0}=0$ and
\begin{equation}\label{eq1.51}
\lambda_n =\Pi^\tau_b(\sigma
0^\infty)=\sum_{j=1}^{k}\tau(\sigma|_j)b^{j-1}+\sum_{i=k+1}^{N_n}\tau(\sigma
0^{i-k})b^{i-1}
\end{equation}
 where $\sigma= \sigma_1\cdots\sigma_k$ is the $q-$adic expansion of $n$ in (\ref{eq1.5}), $\tau(\sigma 0^{N_n-k})\ne 0$
  and $\tau(\sigma 0^{n})= 0$ for all $n>N_n-k$. Under this ordering,
 we have $\Lambda(\tau)=\{\lambda_n\}_{n=0}^\infty$. Let
$$
N^{\ast}_n = \#\left\{\ell: \ \sigma = \sigma_1\cdots\sigma_k, \
\tau(\sigma0^\ell)\neq0 \right\},
$$
where we denote $\#A$ the cardinality of the set $A$.
 The growth rate of $N^{\ast}_n$ is crucial in determining whether $r\Lambda(\tau)$ is a spectrum of $L^2(\mu)$. To describe the growth rate, we let
  ${\mathcal N}^*_{m,n} = \max\{N^{\ast}_k:  q^m\leq k< q^{n}\}$,  ${\cal L}^*_n = \min\{N^{\ast}_k:  q^n\leq k< q^{n+1}\}$ and
  ${\cal M}_n=\max \{N_k: 1\le k<q^n\}$. We have the following two criteria depending on the growth rate of $N^{\ast}_n$.

\medskip

\begin{theorem}\label{th1.7}
 Let $\Lambda= r \Lambda(\tau)$  for a regular
 mapping $\tau$. Then we can find $0<c_1<c_2<1$ so that the
following holds.

\noindent{\rm(i)} If there exists  a strictly increasing sequence $\alpha_n$ satisfying
\begin{equation}\label{eq3.4}
\alpha_{n+1}-{\cal M}_{\alpha_n}\rightarrow\infty, \ \mbox{and} \ \sum_{n=1}^{\infty}c_{1  }^{{\mathcal N}^{\ast}_{\alpha_n,\alpha_{n+1}}} = \infty,
\end{equation}
then $\Lambda$ is a spectrum of $L^2(\mu)$.

\medskip

\noindent{\rm(ii).} If $\sum_{n=1}^{\infty}c_{2}^{{\cal L}_n^{\ast}}<\infty$, then $\Lambda$ is  not a spectrum of $L^2(\mu)$.
\end{theorem}

\medskip

The following is the most important example of Theorem \ref{th1.7}.

\begin{Example}\label{rem1.8}{\rm For a regular mapping $\tau$, if $M: = \sup_{n}\{N^{\ast}_n\}$ is finite, then $\Lambda$ must be a spectrum}.
\end{Example}
\medskip
\proof {\rm Note that ${\mathcal N}^{\ast}_{m,n}\leq M$ and
therefore any strictly increasing sequences $\alpha_n$ will satisfy
the second condition in
 (\ref{eq3.4}). Let $\alpha_1=1$ and  $\alpha_{n+1} = n+{\cal M}_{\alpha_n}$ for $n\ge
 1$. Then the first condition holds and hence $\Lambda$ must be a spectrum by Theorem \ref{th1.7}.}
\eproof

\medskip

We will also see that when there is some slow growth in ${\mathcal N}^{\ast}_n$, $\Lambda$ can still be a spectrum (see Example \ref{example3.1}).
 The exact growth rate for $\Lambda$ to be a spectrum is however hard to obtain from the techniques we used.

%\begin{theorem}
%Suppose $\Lambda= r \Lambda(\tau)$ is a maximal orthogonal set of $L^2(\mu)$ for some regular maximal mapping $\tau$ and
%$$
%\sup_{n\in{\Bbb N}}\{N^{\ast}_n\}<\infty,
%$$
%then $\tau$ is a spectral mapping. (i.e. $\Lambda$ is a spectrum of $L^2(\mu)$)
%\end{theorem}
%
%\medskip
%
%Let $L_n = \min\{N^{\ast}_k:  q^n\leq k<q^{n+1}-1\}$
%\begin{theorem}\label{th1.8}
%Suppose $\Lambda= r \Lambda(\tau)$ is a maximal orthogonal set of $L^2(\mu)$ for some regular maximal mapping $\tau$. Then there exists $0<b<1$ such that if
%$$
%\sum_{n=1}^{\infty}b^{L_n}<\infty,
%$$
% we have $\tau$ is  not a spectral mapping. (i.e. $\Lambda$ is a not spectrum of $L^2(\mu)$)
%\end{theorem}

Now, we can construct some spectra which can have zero Beurling
dimension from regular  orthogonal sets using Theorem \ref{th1.7}.
In fact, they can even be arbitrarily sparse.

\begin{theorem}\label{th1.9}
Let $\mu = \mu_{q,b}$ be a measure defined in (\ref{eq1.1}) with
$b>q$ and $\gcd(q, b)=q$. Then given any increasing non-negative
function $g$ on $[0, \infty)$, there exists a spectrum $\Lambda$ of
$L^2(\mu)$ such that
\begin{equation}\label{eq1.6}
\lim_{R\rightarrow\infty}\sup_{x\in{\Bbb R}}\frac{\#(\Lambda\cap(x-R,x+R))}{g(R)}=0.
\end{equation}
\end{theorem}

\medskip

We also make a study on the irregular spectra, although most
interesting cases are from the regular one.  Let $\tau$ be a maximal
mapping such that it is irregular on $\{I_1 0^\infty, \ldots, I_N
0^\infty\}$, where $I_i\in\Sigma^{\ast}$ and the last word in $I_i$
is non-zero, and is regular on the others in $\Gamma_q$. We define the
corresponding {\it regularized mapping} $\tau_R$:
$$
\tau_{R}(\sigma) =\left\{
  \begin{array}{ll}
    0, & \hbox{if $\sigma=I_i0^{k}$ for $k\ge 1$;} \\
    \tau(\sigma), & \hbox{otherwise.} \\
      \end{array}
\right.
$$

 Our result is as follows:

\begin{theorem}\label{th1.10}
Let $\tau$ be an irregular maximal mapping of $\mu$. Suppose $\tau$
is irregular only on finitely many  $\sigma$ in $\Gamma_q$. Then
$\tau$ is a spectral mapping if and only if a corresponding
regularized mapping $\tau_R$ is a spectral mapping.
\end{theorem}

We will prove this theorem more generally in Theorem \ref{th6.1} by showing the spectral property is not affected if we alter only finitely many
 elements in $\Gamma_q$. However, we don't know whether the same holds if the finiteness assumption on irregular elements is removed.

\bigskip

\section{Maximal orthogonal sets}
In this section, we discuss the existence of orthogonal sets for $\mu_{q,b}$, in particular, Theorem \ref{th1.1} and Theorem \ref{th1.6} are proved.

\medskip

\noindent{\bf Proof of Theorem \ref{th1.1}.} Let gcd$(q,b)=d$.
Suppose $q$ and $b$ are relatively
prime i.e. $d=1$. Let
$$
{\mathcal Z}_n: = \left\{\frac{b^{n}}{q}a:q\nmid a\right\}.
$$
It is easy to see that ${\mathcal Z}(\widehat{\mu}_{q,b}) =\bigcup_{n=1}^{\infty}{\mathcal Z}_n $. Note that for any $a$ with $q\nmid a$,
we have $q\nmid ba$ since gcd$(q,b)=1$. Hence, if $n>1$,
$$
\frac{b^{n}}{q}a = \frac{b^{n-1}}{q}(ba)\in{\mathcal Z}_{n-1}.
$$
This implies that ${\mathcal Z}_1\supset{\mathcal Z}_2\supset\cdots$ and ${\mathcal Z}(\widehat{\mu}_{q,b}) = {\mathcal Z}_1$.  Let
$$
Y_i= \left\{\frac{b}{q}a: q\nmid a, \ a\equiv i \ (\mbox{mod} \ q)\right\},
$$
then ${\mathcal Z}(\widehat{\mu}_{q,b}) = \bigcup_{i=1}^{q-1}Y_i$.
If there exists a mutually orthogonal set $\Lambda$ for $\mu_{q,b}$
with $\#\Lambda\ge q$, we may assume $0\in \Lambda$ so that
$\Lambda\setminus\{0\}\subset {\mathcal Z}(\widehat{\mu}_{q,b})$.
Hence there exists $1\leq i\leq q-1$ such that $Y_i\cap \Lambda$
contains more than 1 elements, say $\lambda_1,\lambda_2$. But then
$\lambda_1-\lambda_2 = \frac{b}{q} r$ where $q|r$. This contradicts
the orthogonal property of $\Lambda$.

\medskip

Suppose now $d>1$, we know $d\leq q$. We first consider $d=q$ and prove that  the measure is a spectral measure. This shows also the second statement.
 Write now $b=qr$ and define ${\mathcal D} =\{0,1,\cdots, q-1\}$ and ${\mathcal S} = \{0,r,\cdots (q-1)r\}$.  Then it is easy to see that the matrix
$$
H: = [e^{2\pi i \frac{ijr}{b}}]_{0\leq i,j\leq q-1} = [e^{2\pi i \frac{ij}{q}}]_{0\leq i,j\leq q-1}
$$
is a Hadamard matrix (i.e. $HH^{\ast}=qI$). This shows ${\frac 1b \mathcal D}$ and ${\mathcal S}$ forms a compatible pair as in \cite{[LaW]}.
Therefore it is a spectral measure by Theorem 1.2 in \cite{[LaW]}.

\medskip

Suppose now $1<d<q$. We have shown that $\mu_{d,b}$ is a spectral measure and hence ${\mathcal Z}(\widehat{\mu}_{d,b})$ contains an infinite bi-zero
 sets $\Lambda$ (i.e. $\Lambda-\Lambda\subset {\mathcal Z}(\widehat{\mu}_{d,b})\cup \{0\}$). We claim that ${\mathcal Z}(\widehat{\mu}_{d,b})
 \subset {\mathcal Z}(\widehat{\mu}_{q,b})$ and hence $ {\mathcal Z}(\widehat{\mu}_{q,b})$ has infinitely many orthogonal exponentials.
 To justify the claim, we write $q=dt$. Note that for $d\nmid a$,
$$
\frac{b^{n}}{d}a = \frac{b^{n}}{q}(ta).
$$
As  $q$ cannot divide $ta$. Hence, $\frac{b^{n}}{q}(ta)\in{\mathcal
Z}(\widehat{\mu}_{q,b})$. This also completes
 the proof of Theorem \ref{th1.1}.
\eproof

\medskip

 \begin{Rem}\label{rem2.1}{\rm In view of Theorem \ref{th1.1}, we cannot decide whether there are spectral measures when
 $1<\gcd(q, b)<q$.
 In general, $\mu_{q,b}$ is the convolutions of several self-similar measures with some are spectral and some are not spectral. If $q=4, \ b=6$,
  we know that $\{0,1,2,3\} = \{0,1\}\oplus\{0,2\}$ and hence}
$$
\widehat{\mu}_{4,6}(\xi) = \prod_{j=1}^{\infty}\left(\frac{1+e^{2\pi i 6^{-j}\xi}}{2}\right)\cdot\prod_{j=1}^{\infty}
\left(\frac{1+e^{2\pi i 2 6^{-j}\xi}}{2}\right)=\widehat{\nu}_{1}(\xi)\widehat{\nu}_{2}(\xi)
$$
{\rm  where $\nu_1 = \mu_{2,6}$ and $\nu_2$ is the equal weight
self-similar measure defined by the IFS with maps $\frac{1}{6}x$ and
$\frac{1}{6}(x+2)$. Hence, $\mu_{4,6}=\nu_1\ast\nu_2$. It is known
that both $\nu_1$ and $\nu_2$ are spectral measures, but we don't
know whether $\mu_{4,6}$ is a spectral measure. If $q=6$ and $b=10$,
then $\{0,1,\cdots, 5\} = \{0,1\}\oplus\{0,2,4\}$ and hence
$\mu_{6,10} $ is the convolution of $\mu_{2,10}$ with a non-spectral
measure with $3$ digits and contraction ratio $1/10$. Because of its
convolutional structure, it may be a good testing ground for
studying the {\L}aba-Wang conjecture \cite{[LaW]} and also for
finding non-spectral measures with Fourier frame \cite{[DL],
[HLL]}.}
\end{Rem}
\medskip

\noindent{\bf Proof of Theorem \ref{th1.6}.} Suppose $\Lambda = r\Lambda(\tau)$ for some maximal mapping $\tau$. We show that it is an maximal
 orthogonal set for $L^2(\mu)$. To see this, we first show $\Lambda$ is a bi-zero set. Pick $\lambda,\lambda'\in\Lambda$, by the definition
 of $\Lambda(\tau)$, we can find two distinct $\sigma$,$\sigma'$  in $\Gamma_q$ such that
$$
\lambda = \frac{b}{q}\Pi^\tau_b(\sigma), \ \lambda' = \frac{b}{q}\Pi^\tau_b(\sigma').
$$
Let $k$ be the first index such that $\sigma|_k\neq\sigma'|_k$. Then for some integer $M$, we can write
 $$
q\lambda -q\lambda' = b\sum_{i=k}^{\infty}(\tau(\sigma|_i)-\tau(\sigma'|_i))b^{i-1} = b^k\left((\tau(\sigma|_k)-\tau(\sigma'|_k)) +b M\right).
$$
By (ii) in Definition \ref{def1.5}, $\tau(\sigma|_k)$ and $\tau(\sigma'|_k)$ are in distinct residue class of $q$. This means $q$ does not
divide $\tau(\sigma|_k)-\tau(\sigma'|_k)$. On the other hand, $q$ divides $b$. Hence, $q$ does not divide $(\tau(\sigma|_k)-\tau(\sigma'|_k)) +b M$.
By (\ref{eq1.3}),  $\lambda-\lambda'$ lies in ${\mathcal Z}(\widehat{\mu})$.

\medskip

To establish the maximality of the orthogonal set $\Lambda$, we show by contradiction. Let $\theta\not\in\Lambda$ but $\theta$ is orthogonal to all
elements in $\Lambda$ . Since $0\in\Lambda$, $\theta\neq0$ and $\theta = \theta-0\in {\mathcal Z}(\widehat{\mu})$. Hence, by (\ref{eq1.3}) we may write
$$
\theta = r(b^{k-1}a),
$$
 where $q$ does not divide $a$. Expand $b^{k-1}a$ in $b-$adic expansion using digits $\{-1,0,\cdots,$$b-2\}$
$$
b^{k-1}a = \epsilon_{k-1}b^{k-1}+\epsilon_{k}b^{k}+\cdots+\epsilon_{k+\ell}b^{k+\ell},
$$
 $q$ does not divide $\epsilon_{k-1}$. Note that there exists unique $\sigma_{s}$, $0\le \sigma_s\le q-1$,  such that $\epsilon_{s}\equiv \sigma_s$
 (mod $q$) for $k-1\le s\le k+\ell$. Denote $\sigma_s=\epsilon_s=0$ for $s>k+\ell$. Since $\theta\not\in \Lambda$,  we can find  the smallest
 integer $\alpha$  such that $\tau(0^{k-2}\sigma_{k-1}\sigma_k\cdots \sigma_{\alpha})\ne \epsilon_\alpha$.
 By (iii) in the definition of $\tau$, we can find $\sigma\in\Gamma_q$ so that $\sigma= 0^{k-2}\sigma_{k-1}\sigma_k\cdots
 \sigma_{\alpha}\sigma'0^{\infty}$ and $\tau$ is regular on $\sigma$, then there exists $M'$ such that
$$
\theta-r\Pi^\tau_b(\sigma)=rb^{\alpha}(\epsilon_{\alpha}-\tau(0^{k-1}\sigma_{k-1}\cdots \sigma_\alpha)+M'b).
$$
By (ii) in the definition of $\tau$, $\tau(0^{k-1}\sigma_{k-1}\cdots \sigma_\alpha)\equiv \sigma_{\alpha}$ (mod $q$), which is also congruent
to $\epsilon_{\alpha}$ by our construction. This implies $\theta-r\Pi^\tau_b(\sigma)$ is not in the zero set of $\widehat{\mu}$ since $q$
divides $\epsilon_{\alpha}-\tau(0^{k-1}i_{k-1}\cdots i_\alpha)$ and $b$ does not divide it either. It contradicts to $\theta$ being orthogonal
 to all $\Lambda$.

\bigskip

Conversely, suppose we are given a maximal orthogonal set $\Lambda$ of $L^2(\mu)$ with $0\in\Lambda$. Then $\Lambda\subset {\mathcal Z}(\widehat{\mu})$.
 Hence, we can write
$$
\Lambda = \{ra_{\lambda}: \lambda\in\Lambda,  a_{\lambda} = b^{k-1}m \ \mbox{for some $k\geq1$ and $m\in{\Bbb Z}$ with $q\nmid m$}\},
$$
where $a_0=0$. Now, expand $a_{\lambda}$ in $b-$adic expansion with digits chosen from ${\mathcal C} = \{-1,0,\cdots,b-2\}$.
\begin{equation}\label{eq2.1}
a_{\lambda} = \sum_{j=1}^{\infty}\epsilon_{\lambda}^{(j)}b^{j-1}.
\end{equation}
 Let $D(\vartheta) = \{\epsilon_{\lambda}^{(1)}:\lambda\in\Lambda\}$ be all the first coefficients of $b$-adic expressions \eqref{eq2.1} of elements
 in $\Lambda$,  and let $
D(c_1,\cdots,c_n) = \{\epsilon_{\lambda}^{(n+1)}: \epsilon_{\lambda}^{(1)}=c_1,\cdots \epsilon_{\lambda}^{(n)} = c_n, \ \lambda\in\Lambda \}
$ be all the  $n+1$-st coefficients of elements in $\Lambda$ whose first $n$ coefficients are fixed, where $c_1, c_2, \ldots, c_n\in\C$. We need the
following lemma.
\medskip

\begin{Lem}\label{lem2.1}
With the notations above, then $D(\vartheta)$  contains exactly $q$ elements which are in distinct residue class (mod $q$) and $0\in D(\var)$. Moreover,
if $D(c_1,\cdots,c_n)$ with all $c_i\in\C$ is non-empty, then it contains exactly $q$ elements which are in distinct residue class (mod $q$) also.
In particular, $0\in D(c_1,\cdots,c_n)$ if $c_1=\cdots=c_n=0$ for $n\ge 1$.
\end{Lem}

\proof  Clearly, by \eqref{eq2.1}, $0\in D(\vartheta)$ and $0\in
D(c_1,\cdots,c_n)$ if $c_1=\cdots=c_n=0$ for $n\ge 1$. Suppose
the number of elements in $D(\vartheta)$ is strictly
less than $q$. We let $\alpha\in \C\setminus D(\vartheta)$ such that
$\alpha$ is not congruent to any elements in $D(\vartheta)$. Then,
for any $\lambda\in \Lambda$, by \eqref{eq2.1} we have
$$
r\alpha-\lambda = r\left(\alpha-\sum_{j=1}^{\infty}\epsilon_{\lambda}^{(j)}b^{j-1}\right)=\frac bq \left(\alpha-\epsilon_\lambda^{(1)}+
\sum_{j=2}^{\infty}\epsilon_{\lambda}^{(j)}b^{j-1}\right).
$$
 Note that $q\nmid (\alpha-\epsilon_\lambda^{(1)})$ for all $\lambda\in \Lambda$ by the assumption, this implies $r\alpha$ is mutually orthogonal to
  $\Lambda$ but is not in $\Lambda$, which contradicts to maximal orthogonality. Hence $D(\vartheta)$ contains at least $q$ elements. If $D(\vartheta)$
  contains more than $q$ elements, then there exists  $a_{\lambda_1}=\sum_{j=1}^{\infty}\epsilon_{\lambda_1}^{(j)}b^{j-1},
 a_{\lambda_2}=\sum_{j=1}^{\infty}\epsilon_{\lambda_2}^{(j)}b^{j-1}$ such that $\epsilon_{\lambda_1}^{(1)}\equiv \epsilon_{\lambda_2}^{(1)}$ (mod $q$)
  and $\epsilon_{\lambda_1}^{(1)}\ne \epsilon_{\lambda_2}^{(1)}$. Then
 $r(a_{\lambda_1}-a_{\lambda_2})=\frac bq (\epsilon_{\lambda_1}^{(1)}-\epsilon_{\lambda_2}^{(1)}+bM)$ for some integer $M$. This means
 $r(a_{\lambda_1}-a_{\lambda_2})$ is not a zero of $\widehat{\mu}$. This contradicts to the mutually orthogonality. Hence, $D(\vartheta)$  contains
 exactly $q$ elements which are in distinct residue class (mod $q$).

\medskip

In general, we proceed by induction. Suppose the statement holds up to $n-1$.  If now $D(c_1,\cdots,c_n)$ is non-empty, then $D(c_1,\cdots,c_k)$ is
 also non-empty for all $k\leq n$. we now show that $D(c_1,\cdots,c_n)$ must contain at least $q$ elements. Otherwise, we consider
 $\theta = r(c_1+\cdots+ c_nb^{n-1}+\alpha b^n)$ where $\alpha$ is in $\C\setminus D(c_1,\cdots,c_n)$ and $\alpha$ is not congruent to any elements
 in $D(c_1,\cdots,c_n)$ (mod $q$). If $\lambda\in \Lambda$ and $\lambda= r(c_1+\cdots+ c_kb^{k-1}+c_{k+1}' b^k+\cdots)$ where $c_{k+1}\neq c_{k+1}'$
 and $k\leq n$, then $c_{k+1}, \ c_{k+1}'\in D(c_1,\cdots c_k)$ and hence $\theta$ and $\lambda$ are mutually orthogonal by the induction hypothesis.
 If $\lambda\in \Lambda$ is such that  the first $n$ digit expansion are equal to $\theta$, the same argument as in the proof of $D(\vartheta)$ shows
  $\theta$ will be orthogonal to this $\lambda$. Therefore, $\theta$ will be orthogonal to all elements in $\Lambda$, a contradiction. Also in a similar
  way as the above, $D(c_1,\cdots,c_n)$ contains exactly $q$ elements can be shown.
\eproof

\medskip

Returning to the proof, by convention, we define $\tau(\vartheta)=0$ and on the first level, we define $\tau(\sigma_1)$ to be the unique element
in $D(\vartheta)$ such that it is congruent to $\sigma_1$ (mod $q$). For  $\sigma = \sigma_1\cdots \sigma_n$, we define $\tau(\sigma \sigma_{n+1})$
to be the unique element in $D(\tau(\sigma|_1),\cdots, \tau(\sigma|_n))$ (it is non-empty from the induction process) that is congruent to
$\sigma_{n+1}$ (mod $q$). Then $\tau(0^k)=0$ for $k\ge 1$.

We show that $\tau$ is a maximal mapping corresponding to $\Lambda$.
(i) is satisfied by above. By Lemma \ref{lem2.1}, $\tau$ is
well-defined with (ii) in Definition \ref{def1.5}. Finally, given a
node $\sigma\in{\Sigma}_q^n$, by the construction of the $\tau$  we
can find $\lambda$ whose first $n$ digits in the digit expansion
(\ref{eq2.1}) exactly equals the value of $\tau (\sigma|_k)$ for all
$1\leq k\leq n$. Since the digit expansion of $\lambda$ becomes $0$
eventually, we continue following the digit expansion of $\lambda$
so that  (iii) in the definition is satisfied.

We now show that $\Lambda= r\Lambda(\tau)$. For each ${a_{\lambda}}$ given in (\ref{eq2.1}), Lemma \ref{lem2.1} with the definition of $\tau$ shows
that there exists unique path $\sigma$ such that  $\tau(\sigma|_n) = \epsilon_{\lambda}^{(n)}$ for all $n$. As the sum is finite, this means
$\Lambda\subset r\Lambda(\tau)$. Conversely, if some $\Pi^\tau_b(\sigma)\in r\Lambda(\tau)$ is not in $\Lambda$, then from the previous proof
we know $\Pi^\tau_b(\sigma)$ must be orthogonal to all elements in $\Lambda$. This contradicts to the maximal orthogonality of $\Lambda$.
Thus, $\Lambda= r\Lambda(\tau)$.

\bigskip

\section{Regular spectra}

 In the rest of the paper, we study under what conditions a maximal orthogonal set is a spectrum or not a spectrum.

\begin{Lem}
\label{lem3.1} \cite{[JP]} Let $\mu$ be a Borel probability
measure in ${\Bbb R}$ with compact support. Then a countable set $\Lambda$
is a spectrum for $L^2(\mu)$ if and only if
$$
Q(\xi):=\sum_{\lambda\in\Lambda}|\widehat{\mu}(\xi+\lambda)|^2\equiv 1,
\quad \mbox{for}\,\, \xi\in{\Bbb R}.
$$
Moreover, if $\Lambda$ is a bi-zero set, then $Q$ is an entire function.
\end{Lem}
Note that the first part of Lemma \ref{lem3.1} is well-known. For the entire function property, we just note that the partial sum
 $\sum_{|\lambda|\leq n}|\cdots|^2$ is an entire function and it is locally uniformly bounded by applying the Bessel's inequality, hence $Q$ is
 entire by the Montel's theorem in complex analysis. One may refer to \cite{[JP]} for the details of the proof.

Let $\delta_a$ be the Dirac measure with center $a$. We define
 $$
\delta_{{\mathcal E}}=\frac1{\#{\mathcal E}}\sum_{e\in{\mathcal E}}\delta_e
$$
for any finite set ${\mathcal E}$. Let
$\mu$ be the self-similar measure in (\ref{eq1.1}). Write $\D=\{0, 1, \ldots, q-1\}$ and $D_k=\frac 1b\D+\cdots+\frac 1{b^k}\D$ for $k\ge 1$.
 We recall that the {\it mask function} of $\D$ is
$$m(\xi)=\frac 1q (1+e^{2\pi i \xi}+\cdots+e^{2\pi i (q-1)\xi})$$
and define $ \mu_{k}=\delta_{{\mathcal D}_k}, $
 then
$$
\widehat{\mu_{k}}(\xi)=\prod_{j=1}^{k}m(b^{-k}\xi)
$$
and it is well-known that $\mu_{k}$  converges weakly to $\mu$ when $k$ tends to
infinity and we have
\begin{equation}\label{eq3.1-}
\widehat{\mu}(\xi)=\widehat{\mu_{k}}(\xi)\widehat{\mu}(\frac{\xi}{b^k}).
\end{equation}

\medskip

\begin{Lem}\label{lem3.2}
Let $\tau$ be a regular  mapping and let $\Lambda =
r\Lambda(\tau)=\{\lambda_k\}_{k=0}^\infty$ be the maximal orthogonal
set determined by $\tau$.  Then for all $n\geq1$,
\begin{equation}\label{eq3.1}
\sum_{k=0}^{q^n-1}\left|\widehat{\mu_n}(\xi+r\lambda_{k})\right|^2\equiv1.
\end{equation}
\end{Lem}

\medskip

\noindent{\bf Proof.}
 We claim that $\{r\lambda_k\}_{k=0}^{q^{n}-1}=\frac{b}{q}\{\lambda_k\}_{k=0}^{q^{n}-1}$ is a spectrum of $L^2(\mu_{n})$. We can then use Lemma \ref{lem3.1} to conclude our lemma.
  Since this set has exactly $q^{n}$ elements, we just need to show the mutually orthogonality. To see this, we note that
\begin{equation}\label{eq3.2}
\widehat{\mu}_n(\xi) = m(\frac{\xi}{b})\cdots m(\frac{\xi}{b^n}).
\end{equation}
Given $l\neq l'$ in $\{0,\cdots q^n-1\}$, let $\sigma
=\sigma_1\cdots \sigma_n$ and $\sigma'=\sigma_1'\cdots \sigma_n'$ be
the $q$-adic expansions of $l$ and $l'$ respectively as in
\eqref{eq1.5}, where $\sigma_n$ and $\sigma_n'$ may be zero. We let
$s\le n$ be the first index such that $\sigma_s\neq \sigma_s'$. Then
we can write
$$
\lambda_{l}-\lambda_{{l'}} =
b^{s-1}(\tau(\sigma|_s)-\tau(\sigma'|_s) +bM)
$$
for some integer $M$. We then have from the integral periodicity of $e^{2\pi i x}$ that
$$
m\left(\frac{r(\lambda_{l}-\lambda_{{l'}})}{b^s}\right)
=m\left(\frac{\tau(\sigma|_s)-\tau(\sigma'|_s)}{q}\right) =0.
$$
It is equal to $0$ because (ii) in the definition of maximal mapping
implies that $q$ does not divide $\tau(\sigma|_s)-\tau(\sigma'|_s)$.
Hence, by (\ref{eq3.2}),
$\widehat{\mu}_n\left(r(\lambda_{l}-\lambda_{{l'}})\right) =0.$
\eproof

\bigskip

Now, we let
$$
Q_n(\xi) = \sum_{k=0}^{q^n-1}\left|\widehat{\mu}(\xi+r\lambda_{k})\right|^2, \mbox{and} \ Q(\xi) = \sum_{k=0}^{\infty}\left|\widehat{\mu}(\xi+
r\lambda_{k})\right|^2.
$$

For any $n$ and $p$, we have the following identity:
\begin{equation}\label{eq4.1}
\begin{aligned}
Q_{n+p}(\xi) =& Q_n(\xi)+ \sum_{k=q^{n}}^{q^{n+p}-1}\left|\widehat{\mu}(\xi+r\lambda_{k})\right|^2\\
             =&Q_n(\xi)+ \sum_{k=q^{n}}^{q^{n+p}-1}\left|\widehat{\mu_{n+p}}(\xi+r\lambda_{k})\right|^2\left|\widehat{\mu}(\frac{\xi+
             r\lambda_{k}}{b^{n+p}})\right|^2.
\end{aligned}
\end{equation}

Our goal is see whether $Q(\xi)\equiv 1$. Then by invoking Lemma
\ref{lem3.1}, we can determine whether we have a spectrum. As $Q$ is
 an entire function by Lemma \ref{lem3.1}, we just need to see the value of $Q(\xi)$ for some small values of $\xi$.  To do this, we need to make a fine
  estimation of the terms  $\left|\widehat{\mu}(\frac{\xi+r\lambda_{k}}{b^{n+p}})\right|^2$ in the above.
Write
$$
c_{\min}=\min\left\{\prod_{j=0}^\infty|m(b^{-j}\xi)|^2: |\xi|\le \frac{b-1}{qb}\right\}>0,
$$
where $|m(\xi)|=\frac{|\sin \pi q\xi|}{q|\sin\pi\xi|}$ and $\prod_{j=0}^\infty|m(b^{-j}\xi)|^2=|m(\xi)\widehat{\mu}(\xi)|^2$. Denote
$$
c_{\max}=\max\left\{|m(\xi)|^2: \frac 1{b^2} \le |\xi|\le \frac{b-1}{qb}\right\}<1.
$$

The following proposition roughly says that the magnitude of the
Fourier transform  is controlled by the number of non-zero digits in
the $b-$adic expansion in a uniform way. Recall that $b=qr$ with
$q,r\ge 2$.

\begin{Prop}\label{prop3.3}
Let $|\xi| \le \frac{r(b-2)}{b-1}$ and let
$$
t = \xi +\sum_{k=1}^{N}d_ib^{n_k},
$$
where $d_i\in\{1,2,\cdots r-1\}$ and $1\leq n_1<\cdots< n_N$. Then
\begin{equation}\label{eq3.3}
c_{\min}^{N+1}\leq\left|\widehat{\mu}(t)\right|^2\leq c_{\max}^N.
\end{equation}
\end{Prop}

\proof First it is easy to check that, for $|\xi|\le \frac{r(b-2)}{b-1}$ and all $d_k\in\{0, 1, 2, \ldots, r-1\}$, we have
\begin{eqnarray}\label{eqxi}
\left|\frac{\xi+\sum_{k=1}^n d_kb^k}{b^{n+1}}\right|&\le& \frac 1{b^{n+1}} \left(\frac{r(b-2)}{b-1}+(r-1)(b+b^2+\cdots+b^n)\right) \nonumber\\
&=&\frac{r(b-2)+(r-1)(b^{n+1}-b)}{b^{n+1}(b-1)} \nonumber\\
&\le&\frac{b-1}{qb}
\end{eqnarray}
for $n\ge1$. The inequality in the last line follows from a direct comparison of the difference and $q\geq 2$.
To simplify notations, we let  $n_0=0$ and $n_{N+1} = \infty$. Then $|\widehat{\mu}(t)|^2$ equals

\begin{equation}\label{eq5.2}
 \prod_{j=1}^{\infty}\left|m\left(b^{-j}t\right)\right|^2
                 =                     \prod_{i=0}^{N}\prod_{j=n_i+1}^{n_{i+1}}\left|m\left(b^{-j}t\right)\right|^2.
\end{equation}
We now estimate the products one by one. By \eqref{eqxi}, we have
$$
\left|\frac{\xi+\sum_{k=1}^{i}d_kb^{n_k}}{b^{n_i+1}}\right|\le \frac{b-1}{qb}.
$$
Hence, together with the integral periodicity of $m$ and the definition of $c_{\min}$, we have for all $i>0$,
\begin{eqnarray}\label{eq5.3}
\prod_{j=n_i+1}^{n_{i+1}}\left|m(b^{-j}t)\right|^2&=&\prod_{j=n_i+1}^{n_{i+1}}\left|m\left(b^{-j}(\xi+
\sum_{k=1}^{i}d_kb^{n_k})\right)\right|^2 \nonumber\\
&\geq&\prod_{j=0}^{\infty}\left|m\left(b^{-j}\left(\frac{\xi+\sum_{k=1}^{i}d_kb^{n_k}}{b^{n_i+1}}\right)\right)\right|^2\geq c_{\min}.
\end{eqnarray}
For the case $i=0$, it is easy to see that
$\left|\frac{\xi}{b}\right|\leq\frac{b-2}{q(b-1)}<\frac{b-1}{qb}.$ Hence, $\prod_{j=n_0+1}^{n_{1}}\left|m(b^{-j}t)\right|^2$
$\geq \prod_{j=0}^{\infty}\left|m\left(b^{-j}(\xi/b)\right)\right|^2\geq c_{\min}$.  Putting this fact and \eqref{eq5.3} into \eqref{eq5.2},
 we have $|\widehat{\mu}(t)|^2\geq c_{\min}^{N+1}$.

\bigskip

We next prove the upper bound. From $|m(\xi)|\leq 1$, \eqref{eq5.2} and the integral periodicity of $m$,
\begin{equation}\label{eq4.8}
|\widehat{\mu}(t)|^2\leq \prod_{i=1}^{N}\left|m\left(b^{-{(n_i+1)}}t\right)\right|^2=\prod_{i=1}^{N} \left|m\left(b^{-{(n_i+1)}}
(\xi +\sum_{k=1}^{i}d_kb^{n_k})\right)\right|^2.
\end{equation}
By \eqref{eqxi} we have
$$|\xi+\sum_{k=1}^id_kb^{n_k}|\ge b^{n_i}-|\xi+\sum_{k=1}^{i-1}d_kb^{n_k}|\ge b^{n_i}-\frac {b^{n_{i-1}}(b-1)}{q}\ge b^{n_i-1}.$$
By \eqref{eqxi}, \eqref{eq4.8},
 the above and the definition of $c_{\max}$, we obtain that $|\widehat{\mu}(t)|^2\le c_{\max}^N$.\eproof

\medskip

 We now prove Theorem \ref{th1.7}. Write $c_1= c_{\min}$ and $c_2=c_{\max}$, where $c_{\min}$ and $c_{\max}$ are in Proposition \ref{prop3.3}.
  Also recall the quantities defined in Section 2. For any $n\in \N$, the $q$-adic expression of $n$ is $\sum_{j=1}^k\sigma_jq^{j-1}$
  with $\sigma_k\ne 0$. Then for the map $\tau$ we have
$$\lambda_n=\sum_{j=1}^k\tau(\sigma_1\cdots \sigma_j)b^{j-1}+\sum_{j=k+1}^{N_n}\tau(\sigma_1\cdots\sigma_k0^{j-k})b^{j-1}$$
where $\tau(\sigma_1\cdots\sigma_k0^{N_n-k})\ne 0$ and
$N^*_n=\#\{\tau(\sigma_1\cdots\sigma_k0^j)\ne 0: k+1\le j\le N_n\}$.
Moreover,
 ${\mathcal N}^{\ast}_{m,n} = \max_{q^m\leq k<q^n}\{N_k^{\ast}\}$, ${\mathcal L}_n^{\ast} = \min_{q^n\leq k<q^{n+1}}\{N_k^{\ast}\}$ and
  ${\mathcal M}_n = \max_{1\leq k<q^{n}}\{N_k\}$.

\medskip

\noindent{\bf Proof of Theorem \ref{th1.7}.}  (i) Let $\alpha_n$ be the increasing sequence satisfying (\ref{eq3.4}) and let $|\xi|\le\frac{b-2}{b-1}$.
Recall (\ref{eq4.1}),
$$
Q_{\alpha_{n+1}}(q^{-1}\xi) =Q_{\alpha_n}(q^{-1}\xi)+ \sum_{k=q^{\alpha_n}}^{q^{\alpha_{n+1}}-1}\left|\widehat{\mu_{\alpha_{n+1}}}(q^{-1}\xi+
r\lambda_{k})\right|^2\left|\widehat{\mu}(\frac{q^{-1}\xi+r\lambda_{k}}{b^{\alpha_{n+1}}})\right|^2.
$$
For $k = q^{\alpha_n},\cdots, q^{\alpha_{n+1}}-1$, we may write $\lambda_k$  as
$$
\lambda_k = \sum_{j=0}^{\alpha_{n+1}-1}c_jb^{j}+ \sum_{j=1}^{M_k}d_jqb^{n_{j}},
$$
where $c_j\in\{-1,\cdots,b-2\}$, $d_j\in\{1,\cdots, r-1\}$ and
$\alpha_{n+1}\leq n_1<n_2<\cdots<n_{M_k}$ with $n_{M_k}=N_k$ and
$M_k\leq N_{k}^{\ast}$, where $N_k$, $N_k^{\ast}$ was defined in
(\ref{eq1.51}) or see the above. Note also that the second term on the
right hand of the above is zero whenever $N_k<\alpha_{n+1}$. Now,
$$
\frac{q^{-1}\xi+r\lambda_{k}}{b^{\alpha_{n+1}}} = \frac{q^{-1}\xi+q^{-1}\sum_{j=1}^{\alpha_{n+1}}c_jb^{j}}{b^{\alpha_{n+1}}}
+\sum_{j=1}^{M_k}d_jb^{n_{j}-\alpha_{n+1}+1}.
$$
Note that
\begin{eqnarray*}\left|\frac \xi q+\frac 1q\sum_{j=1}^kc_jb^{j}\right| \le \frac{b-2}{q(b-1)}+(b-2)\frac{b^{k+1}-b}{q(b-1)} \le \frac{b-2}{q(b-1)}
 b^{k+1}=\frac{r(b-2)}{b-1}b^k
\end{eqnarray*}
for all $k\ge 1$. Hence,
Proposition \ref{prop3.3} implies that
$$
\left|\widehat{\mu}\left(\frac{q^{-1}\xi+r\lambda_{k}}{b^{\alpha_{n+1}}}\right)\right|^2\geq c_1^{1+M_k}\geq c_1^{1+N^{\ast}_k}\geq c_1^{1+
{\mathcal N}^*_{\alpha_n,\alpha_{n+1}}}
$$
for all $q^{\alpha_n}\le k<q^{\alpha_{n+1}}$. Therefore, together with Lemma \ref{lem3.2},
\begin{equation}\label{eq3.5}
\begin{aligned}
Q_{\alpha_{n+1}}(q^{-1}\xi)\geq& Q_{\alpha_{n}}(q^{-1}\xi)+c_1^{1+{\mathcal N}^*_{\alpha_n,\alpha_{n+1}}}\sum_{k=q^{\alpha_n}}^{q^{\alpha_{n+1}}-1}
\left|\widehat{\mu_{{\alpha_{n+1}}}}
(q^{-1}\xi+r\lambda_{k})\right|^2\\
=&Q_n(q^{-1}\xi)+c_1^{1+{\mathcal N}^*_{\alpha_n,\alpha_{n+1}}}\left(1-\sum_{k=0}^{q^{\alpha_n}-1}\left|\widehat{\mu_{{\alpha_{n+1}}}}(q^{-1}\xi
+r\lambda_{k})\right|^2\right).
\end{aligned}
\end{equation}
From elementary analysis, there exists $\delta$, $0<\delta<1$, such
that $|\widehat{\mu}(\xi)|^2$ is decreasing on $(0,\delta)$ (In
fact, since $\widehat{\mu}(0)=1$ and $|\widehat{\mu}(\xi)|^2$ is
entire in complex plane, there exists $\eta>0$ such that
$|\widehat{\mu}(\xi)|<1$ for all $0<\xi<\eta$. If
$|\widehat{\mu}(\xi)|^2$ is not decreasing on $(0,\delta)$ for any
$\delta>0$, we can find a sequence $\xi_n\rightarrow 0$ such that
$(|\widehat{\mu}|^2)'(\xi_n)=0$ and thus
$(|\widehat{\mu}|^2)'\equiv0$ by the entire function property of $|\widehat{\mu}|^2$, this is impossible). In the proof,
it is also useful to note that
$|\widehat{\mu}(-\xi)|=|\widehat{\mu}(\xi)|$. We now argue by
contradiction. Suppose there exists $\Lambda$ such that Theorem
\ref{th1.7} (i) holds but is not a spectrum, then there exists
$t_0<\min\{\delta, \frac{b-2}{b-1}\}$ such that $Q(q^{-1}t_0)<1$
because $Q$ is entire.  For $0\leq k\leq q^{\alpha_n}-1$, we have
$$
\left|\frac{q^{-1}t_0+r\lambda_k}{b^{\alpha_{n+1}}}\right|\leq \frac{1+rb^{{\mathcal M}_{\alpha_n}}}{b^{\alpha_{n+1}}}: = \beta_n.
$$
By the assumption that $\alpha_{n+1}-{\mathcal M}_{\alpha_n}\rightarrow \infty$, we have for all $n$ large, say $n\geq M$, $\beta_n<\delta$ so that
  $\left|\widehat{\mu}\left(\frac{q^{-1}t_0+r\lambda_k}{b^{\alpha_{n+1}}}\right)\right|^2\ge |\widehat{\mu}(\beta_n)|^2$  and we can find $r<1$ such
  that
$$
|\widehat{\mu}(\beta_n)|^{-2}Q(q^{-1}t_0)\leq r<1,\qquad \mbox{for}\ \, n\ge M
$$
because $\beta_n$ tends to zero when $n$ tends to infinity and $\widehat{\mu}(0)=1$.
According to  $\widehat{\mu}(\xi) = \widehat{\mu}_{\alpha_{n+1}}(\xi)\widehat{\mu}(\xi/b^{\alpha_{n+1}})$, we have
\begin{eqnarray*}
|\widehat{\mu}(q^{-1}t_0+r\lambda_k)|^2&=&\left|\widehat{\mu}_{\alpha_{n+1}}(q^{-1}t_0+r\lambda_k)
\widehat{\mu}(\frac{q^{-1}t_0+r\lambda_k}{b^{\alpha_{n+1}}})\right|^2\\
&\ge&|\widehat{\mu}_{\alpha_{n+1}}(q^{-1}t_0+r\lambda_k)|^2|\widehat{\mu}(\beta_n)|^2\\
&\ge&\frac {Q(q^{-1}t_0)}r |\widehat{\mu}_{\alpha_{n+1}}(q^{-1}t_0+r\lambda_k)|^2
\end{eqnarray*}

From (\ref{eq3.5}) and for all $n\geq M$,
$$
\begin{aligned}
Q_{\alpha_{n+1}}(q^{-1}t_0)\geq& Q_{\alpha_n}(q^{-1}t_0)+\left(1-\frac r{Q(q^{-1}t_0)}\sum_{k=0}^{q^{\alpha_n}-1}
|\widehat{\mu}(q^{-1}t_0+r\lambda_k)|^2\right)c_1^{1+{\mathcal N}^*_{\alpha_n,\alpha_{n+1}}}\\
\geq& Q_{\alpha_n}(q^{-1}t_0)+(1-r)c_1^{1+{\mathcal N}^*_{\alpha_n,\alpha_{n+1}}}.
\end{aligned}
$$
Taking summation on $n$ from $M$ to $M+p$ where $p>0$ and note that $Q_n(t)\leq 1$ for any $n$ we have
$$
1\geq Q_{\alpha_{M+p+1}}(q^{-1}t_0)\geq Q_{\alpha_M}(q^{-1}t_0)+(1-r)\sum_{n=M}^{M+p}c_1^{1+{\mathcal N}^{\ast}_{\alpha_n,\alpha_{n+1}}}.
$$
As $\sum_{n=M}^{\infty}c_1^{{\mathcal N}^{\ast}_{\alpha_n,\alpha_{n+1}}}=\infty$ by the assumption, the right hand side of the above tends to infinity.
 This is impossible. Hence, $\Lambda$ must be a spectrum.

%Fixing $n$ and letting $p$ to infinity, we have
%$$
%Q(\xi) \geq Q_n(\xi)+a^N(1-Q_n(\xi)).
%$$
%Take $n$ to infinity also, we have $Q(\xi) = 1$. \qquad$\Box$

\bigskip

\noindent (ii).  The proof starts again at (\ref{eq4.1}) with $p=1$, we have
$$
Q_{n+1}(q^{-1}\xi) =Q_n(q^{-1}\xi)+ \sum_{k=q^{n}}^{q^{n+1}-1}\left|\widehat{\mu_{n+1}}(q^{-1}\xi+r\lambda_{k})\right|^2\left|\widehat{\mu}
(\frac{q^{-1}\xi+r\lambda_{k}}{b^{n+1}})\right|^2.
$$
Since $N^{\ast}_k\geq {\cal L}^*_n$ for $q^n\leq k<q^{n+1}$, so Proposition \ref{prop3.3} implies that
$$
Q_{n+1}(q^{-1}\xi) \leq Q_n(q^{-1}\xi)+ c_2^{{\cal L}^*_n} \sum_{k=q^{n}}^{q^{n+1}-1}\left|\widehat{\mu_{n+1}}(q^{-1}\xi+r\lambda_{k})\right|^2.
$$
Using Lemma \ref{lem3.2} and note that $|\widehat{\mu_{n+1}}(\xi)|^2\geq|\widehat{\mu}(\xi)|^2$, we have
$$
\begin{aligned}
Q_{n+1}(q^{-1}\xi) \leq& Q_n(q^{-1}\xi)+ c_2^{{\cal L}^*_n} \left(1- \sum_{k=0}^{q^{n}-1}\left|\widehat{\mu_{n+1}}(q^{-1}\xi+r\lambda_{k})
\right|^2\right)\\
\leq& Q_n(q^{-1}\xi)+c_2^{{\cal L}^*_n}(1-Q_n(q^{-1}\xi)).\\
\end{aligned}
$$
Hence,
\begin{equation}\label{eq4.7}
1-Q_{n+1}(q^{-1}\xi) \geq (1-Q_n(q^{-1}\xi))(1-b^{{\cal L}^*_n})\geq \cdots \geq(1-Q_1(q^{-1}\xi))\prod_{k=1}^{n}(1-c_2^{{\cal L}^*_k})
\end{equation}
Since $\sum_{n}c_2^{{\cal L}^*_n}<\infty$, $B :=\prod_{k=1}^{\infty}(1-c_2^{{\cal L}^*_n})>0$ and  hence as $n$ tends to infinity in (\ref{eq4.7}),
we have
$$
1-Q(q^{-1}\xi)\geq (1-Q_1(q^{-1}\xi))\cdot B>0.
$$
Therefore, $\tau$ is not a spectral mapping.
\qquad$\Box$

\bigskip
As known from Remark \ref{rem1.8}, $\tau$ is a spectral mapping if
$\sup\{N^{\ast}_n\}$ is finite. Now, we give an example of spectrum
with slow growth rate of $N^{\ast}_n$.

\begin{Example} \label{example3.1} Let $\tau$ be a regular
mapping so that  $N_n \le \log_qn+\log_{c_1^{-2}}\log_q n$ and
$N^{\ast}_n \le \log_{c_1^{-2}}\log_q n$ for $n\ge 1$, where $c_1$
is given in Theorem \ref{th1.7}.  Then $r\Lambda(\tau)$ is a
spectrum of $L^2(\mu)$.
\end{Example}

\medskip

\proof Take $\alpha_n = n^2$. Recall  $\M_n=\max_{1\leq k< q^n}N_k$,
we have
$$\alpha_{n+1}-\M_{\alpha_n}\ge (n+1)^2-n^2-\log_{c_1^{-2}} n^2,$$
 which tends to infinity when $n$ tends to infinity. Note that
$${\mathcal N}_{\alpha_n, \alpha_{n+1}}^\ast=\max_{q^{\alpha_n}\le k<q^{\alpha_{n+1}}} {\cal N}_k^\ast\le \log_{c_1^{-2}}\log_q q^{(n+1)^2}
=\log_{c_1^{-2}}(n+1)^2.$$
Then
$$\sum_{n=1}^\infty c_1^{{\mathcal N}_{\alpha_n, \alpha_{n+1}}^\ast}\ge \sum_{n=1}^\infty c_1^{\log_{c_1^{-2}}(n+1)^2}=\sum_{n=1}^\infty\frac1{n+1}=\infty.$$
By Theorem \ref{th1.7} the result follows.\eproof

%
%
% Choose $\alpha_s=q^s$. Then
%$$\alpha_{s+1}-M_{\alpha_s}= q^{s+1}-\max_{1\le n<q^{q^s}}N_n\ge q^{s+1}-c\log_q q^{q^s}=q^s(q-c),$$
%which tends to infinity when $s$ tends to infinity. Note that
%$$N^*_{\alpha_s, \alpha_{s+1}}=\max\{N^*_n: q^{q^s}\le n<q^{q^{s+1}}\}\le\log_{c_1^{-1}}\log_q\log_q q^{q^{s+1}}=\log_{c_1^{-1}}(s+1).$$
%Hence
%$$\sum_{s=1}^\infty c_1^{N^*_{\alpha_s, \alpha_{s+1}}}\ge \sum_{s=1}^\infty c_1^{\log_{c_1^{-1}}(s+1)}=\sum_{s=1}^\infty\frac 1{s+1}=\infty.$$

\medskip

On the other hand, if $N^{\ast}_n$ is so that  ${\cal L}^{\ast}_n
\ge (1+\epsilon)\log_{c_2^{-1}}n$, for some $\epsilon>0$ and $n\ge
1$, then $r\Lambda(\tau)$ is not a spectrum. This is done by
checking the condition of Theorem \ref{th1.7}(ii) using the similar
method in the above. We therefore omit its detail. Finally, we prove
Theorem \ref{th1.9}.

\medskip

\noindent{\bf Proof of Theorem \ref{th1.9}.}  For any $n\in\N$, $n$ can be expressed as
\begin{equation}\label{finial}
n=\sum_{j=1}^k\sigma_jq^{j-1},
\end{equation}
where all $\sigma_j\in\{0, 1, \ldots, q-1\}$ and $\sigma_k\ne 0$.
Let $\{m_n\}_{n=1}^\infty$ be a strictly increasing sequence of
positive integers with $m_1\ge 2$. We now define a maximal mapping
in terms of this sequence by $\tau(\var)=\tau(0^k)=0$ for $k\ge 1$
and for $n$ as in (\ref{finial}),
$$\tau(\sigma)=\left\{
                                 \begin{array}{ll}
                                   \sigma_k, & \hbox{if}\, \sigma =\sigma_1\cdots\sigma_k, \ \sigma_k\neq 0; \\
                                   0, & \hbox{if}\, \sigma = \sigma_1\cdots\sigma_k0^{\ell}, \ \ell\neq m_n; \\
                                  q, & \hbox{if}\, \sigma = \sigma_1\cdots\sigma_k0^{\ell}\,\, \hbox{and}\,\, \ell=m_n.
                                 \end{array}
                               \right.
$$
%and
%$$\tau(\sigma_1\cdots\sigma_k0^\ell)=\left\{
%                                       \begin{array}{ll}
%                                         q, & \hbox{if}\,\, l=m_n-k+1; \\
%                                         0, & \hbox{otherwise.}
%                                       \end{array}
%                                     \right.
%$$
By the definition we have $\lambda_0=0$ and
$$\lambda_n=\sum_{j=1}^{k}\tau(\sigma_1\cdots\sigma_j)b^{j-1}+qb^{m_n},$$
consequently, $N_n^{\ast} =1$ and by Theorem \ref{th1.7}(i) (see
also Remark \ref{rem1.8}), $\Lambda:=\{\lambda_n\}_{n=0}^\infty$ is
a spectrum for $L^2(\mu)$.

\medskip

We now find $\Lambda$ with density in (\ref{eq1.6}) zero by choosing
$m_n$. To do this, we first note that there exists a strictly
increasing continuous function $h(t)$ from $[0, \infty)$ onto itself
such that $h(t)\le g(t)$ for $t\ge 0$ and it is sufficient to
replace $g(t)$  by $h(t)$ in the proof. In this way, the inverse of
$h(t)$ exists, and we denote it by $h^{-1}(t)$.

Now, note that
$$\lambda_n\le q\frac{b^k-1}{b-1}+qb^{m_n}\le (q+1)b^{m_n}.$$
 Hence,
\begin{equation}\label{fin}
\lambda_{n+1}-\lambda_n\ge qb^{m_{n+1}}-(q+1)b^{m_n}\ge b^{m_n+1}.
\end{equation}
 Therefore, we choose $m_n$ so that  $b^{m_n}\ge 2h^{-1}(b^{n+1})$  for all $n\ge 1$. For any $h(R)\ge 1$,
 there exists unique $s\in \N$ such that $b^{s-1}\le h(R)<b^s$. Then
 $$
 \frac{\sup_{x\in\R}\#(\Lambda\cap(x-R, x+R))}{h(R)}\le \frac{\sup_{x\in\R}\#(\Lambda\cap(x-h^{-1}(b^s), x+h^{-1}(b^s)))}{b^{s-1}}.
 $$
 Note  from \eqref{fin} that  the length of the
open intervals $(x-h^{-1}(b^s), x+h^{-1}(b^s))$  is less than
$\lambda_{n+1}-\lambda_n$ whenever $n\ge s$. This implies that the
set $\Lambda\cap(x-h^{-1}(b^s), x+h^{-1}(b^s))$ contains at most one
$\lambda_n$ where $n\ge s$. We therefore have
$$
\sup_{x\in{\Bbb R}}\#(\Lambda\cap(x-h^{-1}(b^s), x+h^{-1}(b^s)))\le s+1.
$$
Thus the result  follows by taking limit.
 \eproof

\bigskip

\section{Irregular spectra}

  Let $\tau$ be a maximal mapping (not necessarily regular) for $\mu = \mu_{q,b}$ with $b=qr$. Given any $I = \sigma_1\cdots\sigma_k\in\Sigma^k_q=\{0, 1, \ldots, q-1\}^k$ with $\sigma_k\neq 0$.
  Define a map $\tau'$  by

$$
\tau'(\sigma) = \left\{
                  \begin{array}{ll}
                    0, & \hbox{$\sigma = I0^{\ell}$ for $\ell\ge 1$;} \\
                    \tau(\sigma), & \hbox{otherwise.}
                  \end{array}
                \right.
$$
 Clearly $\tau'$ is a maximal mapping. The main result is as follows:

\begin{theorem}\label{th6.1}
With the notation above, $\tau$ is a spectral mapping if and only if $\tau'$ is a spectral mapping.
\end{theorem}

This result shows that if we arbitrarily change the value of $\tau$
along an element in $\Gamma_q$ as above, the spectral property of
$\tau$ is
 unaffected. In particular, Theorem \ref{th1.10} follows as a corollary because we can alter the irregular elements one by one using Theorem
 \ref{th6.1} recursively.

\bigskip

We now prove Theorem \ref{th6.1}. Note that we can decompose
\begin{equation}\label{eq5.1}\Gamma_q=\{\sigma 0^\infty: \sigma\in\Sigma^*_q\}=\bigcup_{I\in\Sigma_q^n}I\Gamma_q
\end{equation}
for all $n\ge 1$. And recall that
$$\Lambda(\tau)=\{\Pi^\tau_b(J): J\in\Gamma_q, \tau\,\,\mbox{is regular on}\,\, J\},$$
where $\Pi^\tau_b(J)=\sum_{k=1}^\infty \tau(J|_k)b^{k-1}$. Denote
naturally $\Pi^\tau_b(I)=\sum_{k=1}^n \tau(I|_k)b^{k-1}$ if
$I\in\Sigma_q^n$, and
$\Pi^\tau_{b,\,I}(J)=\sum_{k=1}^\infty\tau(Ij_1\cdots j_k)b^{k-1}$
for $J=j_1j_2\cdots\in\Gamma_q$ where $IJ$ is regular for $\tau$. Define
also
$$\Lambda_I(\tau)=\{\Pi^\tau_{b,\,I}(J): J\in\Gamma_q, \tau\,\, \mbox{is regular on}\,\, J\}.$$
By \eqref{eq5.1} we have
$$\Lambda(\tau)=\bigcup_{I\in\Sigma^n_q}\big(\Pi^\tau_b(I)+b^n\Lambda_I(\tau)\big).$$

The following is a simple lemma which was also observed in \cite{[DHS]}.

\begin{Prop}\label{lemma6.1}
Let $\tau$ be a tree mapping and $n\ge 1$. Then $r\Lambda(\tau)$ is a spectrum for $\mu$ if and only if all $r\Lambda_I(\tau)$, $I\in\Sigma_q^n$,
 are spectra.
\end{Prop}

\proof Recall that $\mu_k$ satisfies $\widehat{\mu}(\xi)=\widehat{\mu}_k(\xi)\widehat{\mu}(b^{-k}\xi)$ and
$\widehat{\mu}_k(\xi)=\prod_{j=1}^k m(b^{-j}\xi)$ where $m(\xi)=\frac 1q\sum_{j=1}^qe^{2\pi i(j-1)\xi}$. Write
 $Q_I(\xi)=\sum_{\lambda\in\Lambda_I(\tau)}|\widehat{\mu}(\xi+r\lambda)|^2$.
Note that
\begin{eqnarray*}
Q(\xi)=\sum_{\lambda\in\Lambda(\tau)}|\widehat{\mu}(\xi+r\lambda
)|^2&=&\sum_{I\in\Sigma_q^n,\,
\lambda\in\Lambda_I(\tau)}|\widehat{\mu}_n(\xi+r\Pi_b^\tau(I)+rb^n\lambda)|^2\left|\widehat{\mu}\left(\frac{\xi+r\Pi_b^\tau(I)}{b^n}
+r\lambda\right)\right|^2\\
&=&\sum_{I\in\Sigma_q^n,\,
\lambda\in\Lambda_I(\tau)}|\widehat{\mu}_n(\xi+r\Pi_b^\tau(I))|^2\left|\widehat{\mu}\left(\frac{\xi+r\Pi_b^\tau(I)}{b^n}
+r\lambda\right)\right|^2\\
&=&\sum_{I\in\Sigma_q^n}|\widehat{\mu}_n(\xi+r\Pi_b^\tau(I))|^2\cdot Q_I\left(\frac{\xi+r\Pi_b^\tau(I)}{b^n}\right).
\end{eqnarray*}
In a similar proof of Lemma \ref{lem3.2}, we have
$$1\equiv \sum_{I\in\Sigma_q^n}|\widehat{\mu}_n(\xi+r\Pi_b^\tau(I))|^2.$$
Invoking Lemma \ref{lem3.1}, the result follows.
\eproof

\medskip

Proposition \ref{lemma6.1} asserted that spectral property are determined by a finite number of nodes. The following two lemmas show that the spectral
 property of a particular node $\sigma$ can be determined by infinitely many of its offsprings and is {\it independent} of the regularity of
 $\sigma0^{\infty}$. These are the key lemmas to the proof of Theorem \ref{th6.1}.

\begin{Lem}\label{lem6.1} Let $I\in\Sigma_q^*$ with $I\ne \vartheta$, the empty word. If $\tau$ is regular on $I0^\infty$, then
$$
\Lambda_I(\tau)=\{\Pi_{b,\,I}^\tau(0^\infty)\}\cup\bigcup_{k=1}^\infty\bigcup_{j=1}^{q-1}\left(\Pi_{b,\,I}^\tau( 0^{k-1}j)+
b^k\Lambda_{I0^{k-1}j}(\tau)\right),
$$
 where $\Pi_{b, \,I}^\tau(0^{k-1}j)=\tau(I0)+\tau(I0^2)b+\cdots+\tau(I0^{k-1}j)b^{k-1}$. If $\tau$ is irregular on $I0^\infty$, then
$$\Lambda_I(\tau)=\bigcup_{k=1}^\infty\bigcup_{j=1}^{q-1}\left(\Pi_{b,\,I}^\tau(0^{k-1}j)+
b^k\Lambda_{I0^{k-1}j}(\tau)\right).$$
\end{Lem}
\proof Check it directly. \eproof

\begin{Lem}\label{lemma6.2}
Let $\tau$ be a maximal mapping and let $I\in\Sigma_q^{\ast}$. Then
$\Lambda_I(\tau)$ is a spectrum of $\mu$ if and only if
$\Lambda_{I0^{k-1}j}(\tau)$  are spectra of $\mu$ for all $k\geq 1$
and $j=1,\cdots, q-1$.
\end{Lem}
\proof The necessity is clear from Proposition \ref{lemma6.1}. We
now prove the sufficiency. Assume  that $\Lambda_{I0^{k-1}j}(\tau)$
are spectra for all $k\geq 1$ and $j=1,\cdots, q-1$. We need to show
that
$Q_I(\xi)=\sum_{\lambda\in\Lambda_I(\tau)}|\widehat{\mu}(\xi+r\lambda)|^2\equiv
1$.

By the integral periodicity of $m$ and
Lemma \ref{lem3.1} which will be used in the  second equality below, we have for all $k\geq 2$,
$$
\begin{aligned}
&\sum_{j=1}^{q-1}|\widehat{\mu}_k(\xi+r\Pi^{\tau}_{b,I}(0^{k-1}j)|^2\\ =&\sum_{j=1}^{q-1}|\widehat{\mu}_{k-1}(\xi+r\Pi^{\tau}_{b,I}
(0^{k-1})|^2\left|\widehat{\mu}_1
\left(\frac{\xi+r\Pi^{\tau}_{b,I}(0^{k-1})+r\tau(0^{k-1}j)b^{k-1}}{b^{k}}\right)\right|^2\\
=&|\widehat{\mu}_{k-1}(\xi+r\Pi^{\tau}_{b,I}(0^{k-1})|^2 \left(1-\left|\widehat{\mu}_1\left(\frac{\xi+r\Pi^{\tau}_{b,I}(0^{k})}{b^{k}}\right)
\right|^2\right)\\
=&|\widehat{\mu}_{k-1}(\xi+r\Pi^{\tau}_{b,I}(0^{k-1})|^2 -|\widehat{\mu}_{k}(\xi+r\Pi^{\tau}_{b,I}(0^{k})|^2.\\
\end{aligned}
$$
If $k=1$, the above becomes
$\sum_{j=1}^{q-1}|\widehat{\mu}_1(\xi+r\Pi^{\tau}_{b,I}(j)|^2= 1
-|\widehat{\mu}_{1}(\xi+r\Pi^{\tau}_{b,I}(0)|^2$.  Now we simplify
the following terms which is corresponding to the unions of the sets
in Lemma \ref{lem6.1},
\begin{eqnarray}
&&\sum_{k=1}^{\infty}\sum_{j=1}^{q-1}\sum_{\lambda\in\Lambda_{\sigma0^{k-1}j}(\tau)}
|\widehat{\mu}(\xi+r\Pi^{\tau}_{b,I}(0^{k-1}j)+rb^k\lambda)|^2\nonumber\\
&=&\sum_{k=1}^{\infty}\sum_{j=1}^{q-1}\sum_{\lambda\in\Lambda_{\sigma0^{k-1}j}(\tau)}
|\widehat{\mu}_k(\xi+r\Pi^{\tau}_{b,I}(0^{k-1}j)|^2\left|\widehat{\mu}\left(\frac{\xi+r\Pi^{\tau}_{b,I}(0^{k-1}j)}{b^k}
+\lambda\right)\right|^2 \nonumber\\
 &=&\sum_{k=1}^{\infty}\sum_{j=1}^{q-1}|\widehat{\mu}_k(\xi+r\Pi^{\tau}_{b,I}(0^{k-1}j)|^2\nonumber\\
 &=&1-\lim_{N\rightarrow\infty}|\widehat{\mu}_{N}(\xi+r\Pi^{\tau}_{b,I}(0^{N})|^2\nonumber\\
 &=&
 1-\prod_{j=1}^{\infty}\left|m\left(\frac{\xi+r\Pi^{\tau}_{b,I}(0^{j})}{b^j}\right)\right|^2.\label{111}
\end{eqnarray}

\medskip
We now divide the proof into two cases.

\noindent{ Case (i).} If $I0^{\infty}$ is regular, then
 $|\widehat{\mu}(\xi+r\Pi^\tau_{b, I}(0^\infty))|^2 = \prod_{j=1}^{\infty}\left|m(\frac{\xi+r\Pi^{\tau}_{b,I}(0^{j})}{b^j})\right|^2$. Hence, by Lemma \ref{lem6.1},
$$
Q_{I}(\xi) = |\widehat{\mu}(\xi+r\Pi^\tau_{b,
I}(0^\infty))|^2+(\ref{111}) \equiv 1.
$$
This shows $r\Lambda_{I}(\tau)$ is a spectrum.

\medskip

\noindent{ Case (ii).} If $I0^{\infty}$ is irregular, then Lemma \ref{lem6.1} shows that
\begin{equation*}
Q_{I}(\xi) = (\ref{111}) =
1-\prod_{j=1}^{\infty}\left|m\left(\frac{\xi+r\Pi^{\tau}_{b,I}(0^{j})}{b^j}\right)\right|^2.
\end{equation*}
Note that, by (ii) in the definition of the maximal mapping, we may
write
\begin{equation}\label{112} r\Pi_{b,I}(0^{n})
=r\sum_{j=1}^n\tau(I0^j)b^{j-1}= r\sum_{j=1}^{n}(s_jq)b^{j-1} =
\sum_{j=1}^{n}s_jb^{j}
\end{equation} for $s_j\in\{0,1,\cdots, r-1\}$. Then
\begin{equation}\label{eq6.2}
Q_{I}(\xi) =
1-\prod_{j=1}^{\infty}\left|m\left(\frac{\xi+r\Pi^{\tau}_{b,I}(0^{j-1})}{b^j}\right)\right|^2.
\end{equation}
 Suppose on the contrary
$Q_{I}(\xi)<1$ for some $\xi>0$. Since $Q_{I}$ is entire, we may
assume $\xi$ is small, say $|\xi|<\frac{r-1}{b-1}$. From
(\ref{eq6.2}), we must have
$\prod_{j=1}^{\infty}\left|m\left(\frac{\xi+r\Pi^{\tau}_{b,I}(0^{j-1})}{b^j}\right)\right|^2>0.$

 For those $n$ such that $s_n\neq 0$ in \eqref{112}.
$$
\frac1b\leq\frac{s_n}{b}\leq\left|\frac{\xi+\sum_{j=1}^{n}s_jb^j}{b^{n+1}}\right|\leq
\frac{r-1}{b(b-1)}<\frac{1}{q(b-1)}.
$$
Hence, letting $c = \max\{|m(\xi)|^2:
\frac1b\leq|\xi|<\frac{1}{q(b-1)}\}<1$, we have
$$
\left|m\left(\frac{\xi+r\Pi^{\tau}_{b,I}(0^{{n}})}{b^{n+1}}\right)\right|^2
=
\left|m\left(\frac{\xi+\sum_{j=1}^{n}s_jb^j}{b^{n+1}}\right)\right|^2\leq
c.
$$
And
$$
\prod_{j=1}^{\infty}\left|m\left(\frac{\xi+r\Pi^{\tau}_{b,I}(0^{j-1})}{b^j}\right)\right|^2=\lim_{N\rightarrow\infty}\prod_{j=1}^{N}\left|m\left(\frac{\xi+r\Pi^{\tau}_{b,I}(0^{j-1})}{b^j}\right)\right|^2\leq\lim_{N\rightarrow\infty}c^{\#\{n:
s_n\neq 0, \ n\leq N\}}.
$$
As $I0^{\infty}$ is irregular, there exists infinitely many
$s_j\neq 0$. The above limit is zero. This is a contradiction and
hence $\Lambda_{I}(\tau)$ must be a spectrum. \eproof

\medskip

\noindent{\bf Proof of Theorem \ref{th6.1}.}
 By the definition of $\tau$ and $\tau'$, $\Lambda_{I'}(\tau) = \Lambda_{I'}(\tau')$
 for all $I'\neq I$ and $I'\in\Sigma_{q}^k$. Moreover, $\Lambda_{I0^{n-1}j}(\tau)=\Lambda_{I0^{k-1}j}(\tau')$
 for all $k\geq 1$ and $j=1,\cdots q-1$. Therefore, if $\tau$ is a spectrum, then $\Lambda_{I'}(\tau')$ are spectra of $\mu$ for all
 $\sigma'\neq \sigma$ and $\sigma'\in\Sigma_q^{k}$ by Proposition \ref{lemma6.1}. On the other hand, $\Lambda_{I0^{k-1}j}(\tau')$ are spectra
 of $\mu$ also as $\tau$ is a spectrum. By Lemma \ref{lemma6.2}, $\Lambda_{I}(\tau')$ is also a spectrum. We therefore conclude
 that $\Lambda(\tau')$ is a spectrum of $\mu$ by Proposition \ref{lemma6.1} again. The converse also holds by reversing the role of $\tau$ and $\tau'$.
 This completes the whole proof.
\eproof
\medskip
\medskip

\end{document}